\documentclass[12pt]{article}
\newcommand{\bu}{\boldsymbol{u}}

\newtheorem{thm}{Theorem}
\usepackage{amsbsy,amsmath,amssymb}
\date{}
\begin{document}

\title{
\protect\vspace{-2cm}\Large\sc
Integrable (3+1)-dimensional systems\\[1mm] with rational Lax pairs}
\author{
{\sc A. Sergyeyev}\\
Silesian University in Opava, Mathematical Institute,\\ Na Rybn\'\i{}\v{c}ku 1, 74601 Opava, Czech Republic\\
E-mail {\tt artur.sergyeyev@math.slu.cz}
}
\maketitle
\begin{abstract}\protect\vspace*{-1cm}
The search for new integrable (3+1)-di\-men\-si\-on\-al partial differential systems is among the most important challenges in the modern integrability theory.
It turns out that 
such a system can be associated
to any pair of rational functions of one variable in general position, as established below using contact Lax pairs introduced in \cite{as}.
\looseness=-1

{\bf Keywords}: {multidimensional integrable systems; dispersionless systems; contact Lax pairs}
\end{abstract}


\protect\vspace*{-0.5cm}
\section{Introduction}\label{int}

Nonlinear systems, in particular integrable ones, undoubtedly play an important role in modern mathematics, mechanics and physics, cf.\ e.g.\
\cite{cd,das,wl,o,as-ro,asro2,as,tw,tww,z} and references therein.
Integrable systems are particularly interesting as they provide an attractive mix of nonlinearity and tractability. As far as integrable partial differential systems are concerned, 
(3+1)-dimensional ones are particularly important, cf.\ e.g.\ \cite{cd,as} and references therein, as we live in four-dimensional spacetime.

Recently, a novel systematic construction of (3+1)-dimensional integrable systems 
was discovered in \cite{as} using a new class of Lax pairs, the contact Lax pairs. 
\looseness=-1

Namely, let $x,y,z,t$ be independent and $u_1,\dots,u_N$ dependent variables. In \cite{as} it was shown that there exist plethora of (3+1)-dimensional integrable systems 
associated with {\em nonlinear contact Lax pairs} of the form
\begin{equation}\label{nLax4D}
\psi_y=\psi_z F(\psi_x/\psi_z,\bu),\quad \psi_t=\psi_z G(\psi_x/\psi_z,\bu),
\end{equation}
where $\psi=\psi(x,y,z,t)$ is a scalar function and $\bu=(u_1,\dots,u_N)^T$. 
The superscript $T$ stands here and below for the transposed matrix while the subscripts $x,y,z,t$ refer to partial derivatives. Throughout the present paper it is tacitly assumed that all functions are sufficiently smooth for all computations to make sense. This approach can be 
made rigorous using the language of differential algebra, cf.\ \cite{as}.\looseness=-1

To any nonlinear contact Lax pair (\ref{nLax4D}) there corresponds \cite{as} a linear nonisospectral contact Lax pair, and vice versa; 
for this reason only nonlinear contact Lax pairs are considered below.\looseness=-1

In \cite{as} two large classes of nonlinear contact Lax pairs 
leading to integrable (3+1)-dimensional systems were found, 
with the functions $F$ and $G$ being 

\noindent 1) polynomials in $p\equiv\psi_x/\psi_z$ of the form
\begin{equation}
\label{polfg}
F=p^{m+1}+\displaystyle\sum\limits_{i=0}^{m} v_i p^i,\quad 
G=p^{n+1}+\displaystyle\frac{n}{m} v_{m} p^{n}
+\displaystyle\sum\limits_{j=0}^{n-1} w_j p^j,
\end{equation}
where $\bu=(v_{0},\dots,v_{m},w_{0},\dots,w_{n-1})^T$; $m$ and $n$
are arbitrary natural numbers (here and below natural numbers refer strictly to positive integers not including zero), so $N=m+n+1$; 

\noindent 2) rational functions of $p$ of the form
\begin{equation}\label{ratfg}
F=\displaystyle\sum\limits_{i=1}^{m}
\frac{a_i}{p-v_i},\quad G=\displaystyle\sum\limits_{j=1}^{n}
\frac{b_j}{p-w_j},
\end{equation}
where $m$ and $n$ again are arbitrary natural numbers, so now $N=2(m+n)$ and
$\bu=(a_1,\dots,a_m,\allowbreak v_1,\dots, v_m,\allowbreak b_1,\dots, b_n,\allowbreak w_1,\dots,w_n)^T$.

In connection with (\ref{ratfg}) one is immediately led to wonder what happens in a just slightly more general case, when $F$ and $G$ are rational functions of $p$ in general position, that is,
\begin{equation}\label{ratgp}
F=a_0+\displaystyle\sum\limits_{i=1}^{m}
\frac{a_i}{p-v_i},\quad G=b_0+\displaystyle\sum\limits_{j=1}^{n}
\frac{b_j}{p-w_j},
\end{equation}
so now $N=2(m+n+1)$ and
$\bu=(a_0,a_1,\dots,a_m,\allowbreak v_1,\dots, v_m,\allowbreak b_0,b_1,\dots, b_n,\allowbreak w_1,\dots,w_n)^T$.

The goal of the present paper is to answer the natural question posed above. Namely, Theorem~\ref{ratfg-thm} below
introduces a change of variables which turns the nonlinear contact Lax pair associated with $F$ and $G$ from (\ref{ratgp}) into the one associated with $F$ and $G$ from (\ref{ratfg}) under essentially a single assumption of compatibility of the former Lax pair. It turns out that $a_0$ and $b_0$ are basically the artifacts of the gauge freedom, and the said change of variables just removes this freedom.

Note that a similar phenomenon occurs \cite{z} in (2+1) dimensions (this corresponds to putting $\bu_z=0$ and $\psi_z=1$ in the notation of present paper),
so Theorem~\ref{ratfg-thm} can be seen as a generalization of the relevant result of Zakharov \cite{z} to (3+1) dimensions.\looseness=-1

\section{Rational Lax Pairs}\label{mr}


\begin{thm}\label{ratfg-thm}
Let $F$ and $G$ in (\ref{nLax4D}) be rational functions of $\psi_x/\psi_z$ in general position (\ref{ratgp}),
i.e., (\ref{nLax4D}) has the form
\begin{equation}\label{ralaxgen-thm}
\begin{array}{l}
\psi_y=a_0\psi_z+\psi_z^2\displaystyle\sum\limits_{i=1}^{m}
\frac{a_i}{\psi_x-v_i\psi_z},\quad
\psi_t=b_0\psi_z+\psi_z^2\displaystyle\sum\limits_{j=1}^{n}
\frac{b_j}{\psi_x-w_j\psi_z},
\end{array}
\end{equation}
where $m$ and $n$ are any natural numbers, so we have\\
$\bu=(a_0,\dots,a_m,v_1,\dots,v_m,\allowbreak b_0,\dots,\allowbreak b_n, w_1,\dots,w_n)^T$.

Suppose that (\ref{ralaxgen-thm}) is compatible; then there exists a
`potential' $q$ such that $q_z\not\equiv 0$ and
\begin{equation}\label{q-def}
q_y=a_0 q_z, \quad q_t= b_0 q_z.
\end{equation}

Under the above assumptions, upon the change of variables
\begin{equation}\label{chg}
\begin{array}{l}
\tilde{x}=x,\ \tilde{y}=y,\ \tilde{z}=q,\ \tilde{t}=t,\quad \tilde{\psi}=\psi,
\ \tilde q=z, \\[3mm]
\tilde{a}_i=a_i q_z^2,\ \tilde{v}_i=v_i-\displaystyle\frac{q_x}{q_z},\ i=1,\dots,m,
\\[5mm]
\tilde{b}_j=b_j q_z^2,\ \tilde{w}_j=w_j-\displaystyle\frac{q_x}{q_z},\ j=1,\dots,n,
\end{array}
\end{equation}
the system (\ref{ralaxgen-thm}) after omitting the tildas takes the
form
\begin{equation}\label{ralaxgen-thm-tr}
\psi_y=\psi_z^2\displaystyle\sum\limits_{i=1}^{m}
\frac{a_i}{\psi_x-v_i\psi_z},\quad 
\psi_t=\psi_z^2\displaystyle\sum\limits_{j=1}^{n}
\frac{b_j}{\psi_x-w_j\psi_z},
\end{equation}
i.e., it is nothing but
the nonlinear contact Lax pair associated with $F$ and $G$ from (\ref{ratfg}).
\end{thm}
Before proceeding to the proof note that the compatibility condition for (\ref{ralaxgen-thm-tr}) is \cite{as}
the following (3+1)-dimensi\-on\-al integrable system of $2(m+n)$ equations for $2(m+n)$ unknown functions,
\begin{equation}
\label{rls-thm}
\begin{array}{l}
(v_i)_t-\displaystyle\sum\limits_{l=1}^n\biggl\{
\left(\displaystyle\frac{b_l
v_i}{w_l-v_i}\right)_z 
-\left(\displaystyle\frac{b_l}{w_l-v_i}\right)_x-\displaystyle\frac{2 b_l (v_i)_z}{w_l-v_i}\biggr\}=0,
\\[7mm]
(w_j)_y-\displaystyle\sum\limits_{k=1}^m\biggl\{\left(\displaystyle\frac{a_k}{w_j-v_k}\right)_x
-\left(\displaystyle\frac{a_k
v_j}{w_j-v_k}\right)_z+\displaystyle\frac{2 a_k (w_j)_z}{w_j-v_k}\biggr\}=0,
\\[7mm]
(a_i)_t-\displaystyle\sum\limits_{l=1}^n\biggl\{
\displaystyle\displaystyle\frac{3 a_i
(b_l)_z}{w_l-v_i}
-\displaystyle\frac{3 a_i b_l (w_l)_z}{(w_l-v_i)^2}
+\left(\displaystyle\frac{a_i b_l (2
v_i-w_l)}{(w_l-v_i)^2}\right)_z
-\left(\displaystyle\frac{a_i b_l}{(w_l-v_i)^2}\right)_x\biggr\}=0,\\[7mm]
(b_j)_y-\displaystyle\sum\limits_{k=1}^m\biggl\{\displaystyle\displaystyle\frac{3 a_k
(b_j)_z}{w_j-v_k}
-\displaystyle\frac{3 a_k b_j (w_j)_z}{(w_j-v_k)^2}
+
\left(\displaystyle\frac{a_k b_j (2
v_k-w_j)}{(w_j-v_k)^2}\right)_z-
\left(\displaystyle\frac{a_k b_j}{(w_j-v_k)^2}\right)_x\biggr\}=0,
\end{array}
\end{equation}
where $i=1,\dots,m$ and $j=1,\dots,n$.

The above system is obviously determined; in particular, it has as many equations as dependent variables. Moreover, system (\ref{rls-thm}) can \cite{as} be turned into a system of Cauchy--Kowalewski type by a simple change of independent variables, $X=x$, $Y=y-t$, $Z=z$, $T=y+t$. Integrability of (\ref{rls-thm}) in the sense of existence of {\em linear} Lax pair immediately follows from the general results of \cite{as}.

\noindent{\em Proof of Theorem~\ref{ratfg-thm}.} 
The left-hand side of the compatibility condition for (\ref{nLax4D}), that is,
\begin{equation}\label{cc}
(\psi_y)_t-(\psi_t)_y=0,
\end{equation}
where the derivatives are evaluated by virtue of (\ref{nLax4D}), 
obviously is a rational function of $p\equiv\psi_x/\psi_z$.

Bringing this rational function to a common denominator and
then equating to zero the coefficients at the powers of $p$ in the numerator
yields a system of $2n+2m+1$ PDEs for the $2m+2n+2$ unknown functions $u_A$, 
i.e., this system is underdetermined.

Equating to zero the coefficient at the highest power of $p$ of the said numerator yields the equation
\begin{equation}\label{ab}
(a_0)_t-(b_0)_y-b_0 (a_0)_z+a_0 (b_0)_z=0.
\end{equation}

One can easily verify that a general solution of (\ref{ab}) can be written as 
\begin{equation*} 
a_0=q_y/q_z,\quad b_0=q_t/q_z,
\end{equation*}
where $q$ is an arbitrary function of $x,y,z,t$ such that $q_z\not\equiv 0$. The above system is nothing but (\ref{q-def}),
so $q$
introduced in Theorem~\ref{ratfg-thm} is indeed well-defined if (\ref{ralaxgen-thm}) is compatible.\looseness=-1

The presence of an arbitrary function $q$ is a manifestation of the gauge freedom in the system under study, just as in the
(2+1)-dimensional case studied in \cite{z}.

Straighforward but cumbersome computations show that passing to new variables given by (\ref{chg}) turns,
modulo omitting the tildas, (\ref{ralaxgen-thm}) into (\ref{ralaxgen-thm-tr}).
Notice that (\ref{ralaxgen-thm-tr}) does not involve $a_0$ and $b_0$, i.e.,
using the transformation (\ref{chg}) has removed the gauge freedom
associated with $q$. This remark completes the proof. \hfill 
$\Box$
\looseness=-1

\section{Outlook}

It is immediate from the above that there exist non-over\-de\-ter\-mined (3+1)-di\-men\-si\-on\-al integrable systems with Lax pairs (\ref{ralaxgen-thm}) associated with arbitrary pairs of rational functions of a single variable in general position.\looseness=-1

These systems possess 
hidden gauge freedom just like their (2+1)-dimensional counterparts with Lax pairs of the form
\[
\begin{array}{l}
\psi_y=a_0+\displaystyle\sum\limits_{i=1}^{m}
\frac{a_i}{\psi_x-v_i},\quad 
\psi_t=b_0+\displaystyle\sum\limits_{j=1}^{n}
\frac{b_j}{\psi_x-w_j}
\end{array}
\]
studied 
in \cite{z}. Moreover, upon the removal of the said gauge freedom through the change of variables presented in Theorem~\ref{ratfg-thm} the systems in question together with their Lax pairs (\ref{ralaxgen-thm}) become equivalent to somewhat simpler systems (\ref{rls-thm}) with Lax pairs (\ref{ralaxgen-thm-tr}) studied in \cite{as}. While the systems associated with (\ref{ralaxgen-thm}) are underdetermined, this is not the case for systems (\ref{rls-thm}) which are equivalent to systems of Cauchy--Kowalevski type as discussed in Section~\ref{mr}.\looseness=-1

The significance of these results consists inter alia in revealing the 
breadth of the class of (3+1)-di\-men\-si\-on\-al integrable systems in general, as well as of the subclass thereof associated with contact Lax pairs (\ref{nLax4D}).\looseness=-1 

Moreover, the results of the present paper immediately lead to a natural open problem of 
going beyond the rational case and 
finding examples of (3+1)-dimensi\-on\-al systems with Lax pairs (\ref{nLax4D}) whose Lax functions $F$ and $G$ have more involved 
dependence on $\psi_x/\psi_z$. 
The author intends to address this in future work.\looseness=-1

\section*{Acknowledgments}
{This research was
supported in part by the Ministry of Education, Youth and Sports of the
Czech Republic (M\v{S}MT \v{C}R) under RVO funding for I\v{C}47813059, and by the
Grant Agency of the Czech Republic (GA \v{C}R) under grant P201/12/G028.}

The author declares he has no conflict of interests.

\end{document}